\documentclass[11pt,a4paper]{article}
\usepackage{amsmath, amssymb, theorem, latexsym,color}
\usepackage{graphicx}

\newtheorem{theorem}{Theorem}[section]

\newtheorem{proposition}[theorem]{Proposition}

\newtheorem{conjecture}[theorem]{Conjecture}

\newtheorem{remark}[theorem]{Remark}

\allowdisplaybreaks

\def\neweq#1{\begin{equation}\label{#1}}
\def\endeq{\end{equation}}

\def\Om {{\Omega}}
\def\la {{\lambda}}

\newcommand {\ar}{\rightarrow}
\newcommand {\pa}{\partial}

\numberwithin{equation}{section}

\title{ Lower bound for the number of critical points  of minimal spectral  $k$-partitions for $k$ large.}

\author{\Large\sc Bernard Helffer \\ [\medskipamount]
\it Laboratoire de math\'ematiques (UMR 8628), Universit\'e Paris-Sud, \\
\it B\^atiment 425, 91405 Orsay Cedex, France, \\[\medskipamount]
\it Laboratoire de Math\'ematiques Jean Leray, Universit\'e de Nantes,\\
E-mail: \texttt{bernard.helffer@math.u-psud.fr}
}

 \date{26 March 2015}

\begin{document}
\maketitle

 \begin{abstract}
 In a recent paper with Thomas Hoffmann-Ostenhof, we proved  that the number of critical points $\nu_k$  in the boundary set of a $k$-minimal partition
  tends to $+\infty$ as $k\ar +\infty$. In this note, we show that $\nu_k$ increases linearly with $k$ as suggested by a hexagonal conjecture about the asymptotic behavior of the energy of these minimal partitions. As the original proof by Pleijel, 
  this involves Faber-Krahn's inequality and Weyl's formula, but this time, due to the magnetic characterization of the minimal partitions, we have to establish 
   a Weyl's formula for Aharonov-Bohm operator controlled with respect to a $k$-dependent number of poles.
 \end{abstract}

\section{Introduction}
We consider  the Dirichlet  Laplacian in a bounded
regular domain $\Omega\subset \mathbb R^2$.
In \cite{HHOT} we have analyzed the  elations between the nodal domains
 of the real-valued  eigenfunctions of this Laplacian and the partitions of
 $\Omega$ 
 by $k$ disjoint 
 open sets $ D_i$ which are minimal in the sense that  the
 maximum over the $ 
D_i$'s of the ground state energy (or smallest eigenvalue) of  the Dirichlet
 realization
 of the Laplacian  in $ D_i$ is minimal. 
  We denote by $ (\lambda_j(\Omega))_{j\in \mathbb N}$ the
non decreasing sequence of its eigenvalues and by $\phi_j$ some associated
 orthonormal basis of real-valued eigenfunctions. 
The groundstate  $ \phi_1$ can be chosen to be
strictly positive in $ \Om$, but the other eigenfunctions 
$  \phi_j$ ($j>1$) must have non empty zeroset in $\Omega$.
By the zero-set of a real-valued continuous function $u$ on 
$ \overline\Om $, we mean
$
N(u)=\overline{\{x\in \Om\:\big|\: u(x)=0\}}
$
and call the components of $ \Om\setminus N(u)$ the nodal
domains of $u$.
The number of 
nodal domains of $ u$ is called $ \mu(u)$. These
$\mu(u)$ nodal
domains define a $k$-partition of $ \Omega$, with $k=\mu(u)$. 

We  recall that the  Courant nodal
 Theorem \cite{Cou} says that, for  $k\geq 1$,  and if  $ E(\lambda_k)$ denotes the eigenspace associated with $\lambda_k\,$, then,   for all real-valued $ u\in  E(\lambda_k)\setminus \{0\}\;,\;
\mu (u)\le k\,.
$

 A  theorem due to Pleijel \cite{Pl} in
1956  says
that this cannot be true when the dimension (here we consider the
$2D$-case) is larger than one. In the next section, we describe the link of these results with the question of spectral minimal partitions which were introduced by Helffer--Hoffmann-Ostenhof--Terracini \cite{HHOT}.

\section{Minimal spectral partitions}
We now introduce for $k\in \mathbb N$ ($k\geq 1$),
 the notion of $k$-partition. We
call {\bf  $ k$-partition}  of $ \Omega$ a family 
$ \mathcal D=\{D_i\}_{i=1}^k$ of mutually disjoint sets in $\Omega$.
 We denote by $ \mathfrak O_k(\Omega)$ the set of open connected
partitions of $\Omega$. 
We now introduce the notion of energy of the partition  $\mathcal D$ by
\begin{equation}\label{LaD}
\Lambda(\mathcal D)=\max_{i}\la(D_i)\,.
\end{equation} 
Then we define for any $k$ the minimal energy in $\Omega$ by
\begin{equation}\label{frakL} 
\mathfrak L_{k}(\Omega)=\inf_{\mathcal D\in \mathfrak O_k}\:\Lambda(\mathcal D).
\end{equation}
and  call  $ \mathcal D\in \mathfrak O_k$ a minimal $k$-partition if 
  $ \mathfrak L_{k}=\Lambda(\mathcal D)$. 
We  associate with a partition its {\bf boundary set}:
\begin{equation}\label{assclset} 
N(\mathcal D) = \overline{ \cup_i \left( \partial D_i \cap \Omega
  \right)}\;.
\end{equation}
The properties of the boundary of a minimal partition are quite close to the properties of nodal sets can be described in the following way: 
\begin{enumerate}
\item
Except for finitely many distinct $ X_i\in \Om\cap N$
 in the neighborhood of which $ N$ is the union of $\nu_i= \nu(X_i)$
smooth curves ($ \nu_i\geq 3$) with one end at $ X_i$,  
$ N$ is locally diffeomorphic to a regular 
curve.\\
\item
$ \pa\Om\cap N$ consists of a (possibly empty) finite set
of points $ Y_i$. Moreover  
 $N$ is near $ Y_i$ the union 
of $ \rho_i$ distinct smooth half-curves which hit
$ Y_i$.\\
\item $ N$  has the {\bf equal angle
  meeting
 property}\footnote{ The half curves meet with equal angle at each critical
 point of $ N$ and also at the boundary together with the
 tangent to the boundary.}
\end{enumerate}
The $X_i$ are called the critical points and define the set
$X(N)$.  A particular role is played by $X^{odd}(N)$ corresponding to the critical points for which $\nu_i$ is odd.  
  
 It has been proved by Conti-Terracini-Verzini (existence)
 and  Helffer--Hoffmann-Ostenhof--Terracini (regularity)  (see \cite{HHOT} and references therein)  that for any $ k$, there exists a  minimal  regular $
k$-partition, and moreover that  
any  minimal $ k$-partition has a  regular
representative\footnote{possibly 
after a modification of the open sets of the partition by capacity $0$ subsets.}.

In a recent paper with Thomas Hoffmann-Ostenhof \cite{HHO},  we proved  that the number of odd critical points  of a minimal  $k$-partition  $\mathcal D_k$
\begin{equation}\label{weak}
 \nu_k:=\# X^{odd} (N(\mathcal D_k))
 \end{equation}
  tends to $+\infty$ as $k \ar + \infty$. \\

 In this note, we will  show that it increases linearly with $k$ as suggested by the hexagonal conjecture as discussed in \cite{BHV,CL, BBO,HHO}. This conjecture says that
\begin{equation}\label{hexa}
A(\Omega)  \lim_{k\ar +\infty} \frac{\mathfrak L_k(\Omega)}{k} = \lambda({\rm Hexa_1})\,,
\end{equation}
where ${\rm Hexa_1}$ denotes the regular hexagon of area $1$ and $A(\Omega)$ denotes the area of $\Omega$.\\
Behind this conjecture, there is the idea that $k$-minimal partitions will look (except at the boundary where one can imagine that pentagons will appear) as the intersection with $\Omega$  of a tiling by hexagons of area $\frac{1}{k} A(\Omega)$.\\
The proof presented here gives not only a better result but is at the end simpler, although based on the deep  magnetic characterization of minimal partitions of \cite{HHmag} which will be recalled in the next section.
 \section{ Aharonov-Bohm operators and magnetic characterization.}
Let us recall some definitions about the Aharonov-Bohm 
Hamiltonian in an open set $\Omega$ (for short $ {{\bf A}{\bf B}}X$-Hamiltonian) with a singularity at $ X\in \Omega$ as considered in \cite{HHOO,AFT}. We denote by $ X=(x_{0},y_{0})$ the coordinates of the pole and 
consider the magnetic potential with  flux at $ X$: 
$ \Phi = \pi  $, defined in $\dot{\Omega_X}= \Omega \setminus \{X\}$:
\begin{equation}
{{\bf A}^X}(x,y) = (A_1^X(x,y),A_2^X(x,y))=\frac 12\, \left( -\frac{y-y_{0}}{r^2}, \frac{x-x_{0}}{r^2}\right)\,.
\end{equation}
The ${{\bf A}{\bf B}}X$-Hamiltonian is defined  by considering the Friedrichs
extension starting from $ C_0^\infty(\dot \Omega_{X})$
 and the associated differential operator is
\begin{equation}
-\Delta_{{\bf A}^X} := (D_x - A_1^X)^2 + (D_y-A_2^X)^2\,\mbox{with }D_x =-i\pa_x\mbox{ and }D_y=-i\pa_y.
\end{equation}
Let $ K_{X}$ be the antilinear operator 
$  K_{X} = e^{i \theta_{X}} \; \Gamma\,$,
with $  (x-x_0)+ i (y-y_0) = \sqrt{|x-x_0|^2+|y-y_0|^2}\, e^{i\theta_{X}}\,$,  $\theta_X$ such that $d\theta_X= 2 {\bf A}^X\,$, 
and where $ \Gamma$ is the complex conjugation operator  
$ \Gamma u = \bar u\,$.  A  function $ u$ is called  $ K_{X}$-real, if 
$ K_{X} u =u\,.$ The operator $ -\Delta_{{\bf A}^X}$ is preserving the 
$ K_{X}$-real functions and we can consider a
basis of $K_{X}$-real eigenfunctions. 
Hence we only  analyze the
restriction  of the ${{\bf A}{\bf B}}X$-Hamiltonian
to the $ K_{X}$-real space $ L^2_{K_{X}}$ where
$$
 L^2_{K_{X}}(\dot{\Omega}_{X})=\{u\in L^2(\dot{\Omega}_{X}) \;,\; K_{X}\,u =u\,\}\,.
$$
This construction can be extended to 
  the case of a
 configuration with $ \ell$ distinct points $ X_1,\dots, X_\ell$ (putting a  flux $ \pi$ at each
 of these points). We just take as magnetic potential 
$$ 
{\bf A}^X = \sum_{j=1}^\ell {\bf A}^{X_j}\,, \mbox{ where } X=(X_1,\dots,X_\ell)\,.$$
We can also construct the antilinear
 operator $ K_X$,  where $ \theta_X$ is replaced by a
 multivalued-function $ \phi_X$ such that $ d\phi_X = 2 {\bf A}^{X}$.  We can then  consider the 
 real subspace of the $ K_X$-real
 functions in $ L^2_{K_{X}}(\dot{\Omega}_{X})$.
 It was shown  in \cite{HHOO} and \cite{AFT} that  the $ K_{X}$-real eigenfunctions have a regular nodal set
 (like the eigenfunctions of the Dirichlet Laplacian) with the
 exception that at each singular point $ X_j$ ($ j=1,\dots,\ell$) 
 an odd number of half-lines  meet.\\
 The next theorem which
  is the most interesting part of the magnetic characterization of the minimal partitions given in \cite{HHmag} will play a basic role in the proof of our main theorem.
\begin{theorem}\label{thchar}[Helffer--Hoffmann-Ostenhof]~\\
Let $ \Omega$ be simply connected. If $ \mathcal D$ is a $k$-minimal partition of $\Omega$, then, by choosing  $ (X_1,\dots, X_\ell)= X^{odd}(N(\mathcal D))$, $ \mathcal D$ is the nodal partition 
of some $ k$-th $ K_{X}$-real eigenfunction  of  the Aharonov-Bohm Laplacian  associated with $ \dot{\Omega}_X$.
\end{theorem}
 
\section{Analysis of the critical sets in the large limit case}  
We can now state our main theorem, which improves \eqref{weak} as proved in  \cite{HHO}.
 \begin{theorem}[Main theorem]~\\
 Let $(\mathcal D_k)_{k\in \mathbb N} $ be a sequence of regular minimal $k$-partitions.
 Then there exists $c_0 >0$ and $k_0$ such that for $k\geq k_0$, 
 $$\nu_k:= \# X^{odd}(N(\mathcal D_k ))\geq c_0 k\,.
 $$
 \end{theorem}
 
 {\bf Proof}\\
 The proof is inspired by the proof of Pleijel's theorem, with the particularity that the operator, which is now the Aharonov-Bohm operator  will depend on $k$.  For each $\mathcal D_k$, we consider the corresponding Aharonov-Bohm operator as constructed in Theorem \ref{thchar}.\\
 We come back to the proof of the lower bound of the Weyl's formula but we will make a partition in squares depending on 
 $$\lambda =\mathfrak L_k\,.$$  
 We introduce a square $Q_p$  of size $t/\sqrt{\lambda}$ with $t\geq 1$ which will be chosen large enough (independently of $k$) and will be determined later. 
 Having in mind the standard proof of the Weyl's formula (see for example \cite{CH}),  we  recall the following proposition
\begin{proposition}\label{Corw}~\\ If $\mathcal D$ is a partition of $\Omega$, then
\begin{equation}\label{mmw}
\sum_i n(\lambda, D_i) \leq n(\lambda, \Omega)\;.
\end{equation}
\end{proposition}
Here $n(\lambda,\Omega)$ is the counting function of the eigenvalues $< \lambda$ of $H(\Omega)$.\\ This proposition  is actually present in the proofs of
 the asymptotics of the counting function.  We will apply this proposition in the case of Aharonov-Bohm operators $H:= -\Delta_{{\bf A}^X }$ restricted to $K_X$-real $L^2$ spaces.  $H(D)$ means the Dirichlet realization (obtained via the Friedrichs extension theorem) of  $-\Delta_{{\bf A}^X}$ in an open set  $D\subset  \Omega$.\\ 
 \begin{remark}\label{remutile}
 Note that if no pole belongs to $D$ (a pole on $\partial D$ is permitted) and if $D$ is simply connected, then $H(D)$ is unitary equivalent  (the magnetic potential can be gauged away) to the Dirichlet Laplacian in $D$.
 We refer to \cite{AFT,Len,NT} for a careful analysis of the  domains of the involved operators.
 \end{remark}
 
We now consider a maximal partition of $\Omega$ with squares $Q_p$ of size  $t/\sqrt{\lambda}$ with the additional rule that the squares should not contain  the odd critical points  of $\mathcal D_k$.\\
The area   $A(\Omega_{k,t,\lambda})$, where   $\Omega_{k,t,\lambda}$ is the union of these squares,   satisfies
  $$
   A( \Omega_{k,t,\lambda} ) \geq  A(\Omega) -  \ell   t^2 
  /\lambda - C (t,\Omega)\frac{1}{\sqrt{k}}\,.
  $$
  The second term on the right hand side estimates from above the area of the squares containing a critical point and the last term takes account of the effect of the boundary.\\
  
 Note that this lower bound of $A(\Omega_{k,t,\lambda})$ leads to the estimate of the cardinal of the squares using 
 $
 \# \{ Q_p \}  =   A( \Omega_{k,t,\lambda} ) \frac{\lambda}{t^2}\,.
 $
 In each of the squares, because (as recalled in Remark \ref{remutile}) the magnetic Laplacian is isospectral to the usual Laplacian, 
 we have (after a dilation argument) :
 $$ n( \lambda , Q_p) = n\left(t, (0,1)^2\right)\,.$$
Hence we need 
  to find a lower bound of $n(t):=n(t, (0,1)^2)$, the number of eigenvalues less than $t^2$ for the standard Dirichlet Laplacian  in the fixed unit square. \\
 We know, that for any $\epsilon >0$ there exists $t$ such that 
 \begin{equation} \label{wq}
 n(t) \geq (1-\epsilon) \frac{1}{4\pi} t^2 \,.
 \end{equation}
 This leads, using Proposition \ref{Corw} for $H= - \Delta_{{\bf A}^X}$ (remember that $X$ is given by the magnetic characterization of $\mathcal D_k$) and applying \eqref{wq} in each square, to the lower bound as $k\ar +\infty$, 
 \begin{equation}\label{z1}
k= n(\mathfrak L_k,\Omega)  \geq \left( \frac{1}{4\pi} (1-\epsilon)t^2\right)\,  \left( A(\Omega) -  \ell  t^2  /\mathfrak L_k  + o(1) \right) \, \left(\frac{\mathfrak L_k}{t^2}\right) \,
\end{equation}
Let us recall from \cite{HHOT} the following consequence of  Faber-Krahn's inequality
\begin{equation}\label{z2}
A(\Omega) \frac{\mathfrak L_k(\Omega)}{k} \geq  \pi {\bf j}^2\,,
\end{equation}
where ${\bf j}\sim 2.405 $ is the first zero of the first Bessel function.

Dividing \eqref{z1}  by $k$ and using \eqref{z2}, 
we get, as $k\ar +\infty$
$$
1 \geq  \frac{{\bf j}^2}{4} (1-\epsilon)  (1 -  \frac{\ell}{k}  t^2 \pi^{-1} {\bf j}^{-2}) (1 + o(1)).
 $$
 If we assume that the number $\ell$  of critical points satisfies
  $$
  \ell \leq \alpha k\,, \, \mbox{for some }\alpha >0\,,
  $$
  we get
 \begin{equation}\label{contra}
1 \geq  \frac{{\bf j}^2}{4} (1-\epsilon)  (1 -  \alpha  t^2 \pi^{-1} {\bf j}^{-2}) (1 + o(1)).
 \end{equation}

 We see that if $\epsilon$ is small enough (this determines $t =t(\epsilon)$) and $\alpha t^2$ is small enough  such that
 $$
\frac{{\bf j}^2}{4} (1-\epsilon)  (1 -  \alpha  t^2 \pi^{-1} {\bf j}^{-2})  >1
 $$
 (this gives the condition on $\alpha$), 
 we will get  a contradiction for $k$ large.\\

As recalled in \cite{HHO},  Euler's formula implies that for a minimal $k$-partition
 $\mathcal D$ of a simply connected domain $\Omega$
 the cardinal of $X^{odd}(N(\mathcal D))$ satisfies
\begin{equation}
\# X^{odd}(N(\mathcal D) ) \leq 2k -4 \,. 
\end{equation}
This estimate seems optimal  and is compatible with the hexagonal conjecture, which, for critical points, will read
\begin{conjecture}
\begin{equation}
\lim_{k\ar +\infty} \frac{\# X^{odd}(N(\mathcal D_k) )}{k}=2 \,. 
\end{equation}
\end{conjecture}

\section{Explicit lower bounds}
 Looking at the proof of the main theorem, the contradiction is obtained if \eqref{contra} is satisfied.
Using the universal lower bound for $n(t)$ (see for example \cite{Pl}), we have, if $t\geq 2$
\begin{equation}
n(t) > \frac {1}{4\pi} t^2 - \frac{2}{\pi^2}  t + \frac{1}{\pi^2}  \,.
\end{equation}
We look for $t=t(\epsilon) \geq 2$ such that
\begin{equation*}
 \frac {1}{4\pi} t^2 - \frac{2}{\pi^2}  t  + \frac{1}{\pi^2} \geq (1-\epsilon) \frac {1}{4\pi} t^2 \,,
\end{equation*}
which leads to the condition
\begin{equation}\label{conde}
\epsilon \frac {1}{4\pi} t^2 - \frac{2}{\pi^2}  t + \frac{1}{\pi^2}\geq 0\,.
\end{equation}
We can choose 
$t(\epsilon) = \max (2, \frac{ 8}{\epsilon \pi })\,$.
We then get a condition on $\alpha$ through \eqref{contra}.   For some admissible $\epsilon$, i.e satisfying: $$ \frac{{\bf j}^2}{4} (1-\epsilon) >1 \,,
$$
the proof works if
$ \alpha < c_0(\epsilon) $, with $c_0(\epsilon) $ solution of 
  \begin{equation} \label{zc0}
  1 = \frac{{\bf j}^2}{4} (1-\epsilon)  (1 -  c_0(\epsilon)  t(\epsilon)^2 \pi^{-1} {\bf j}^{-2}) \,.
  \end{equation}
Hence the $c_0$ announced in the theorem can be chosen as
$$
c_0:= \sup_{\epsilon \in (0, 1-4/{\bf j}^2)}  c_0(\epsilon)\,,
$$
It remains to determine this $\sup$.
Note that $1-4/{\bf j}^2 \sim  0,36\,$.
Hence we can assume $t(\epsilon)= \frac{8}{\epsilon \pi}$ and get for $c_0(\epsilon)$ the equation
\begin{equation}
 c_0(\epsilon)=  \epsilon^2  2^{-6} \pi^3 {\bf j}^{2} \left(1 - \frac{4}{{\bf j}^2(1-\epsilon)}\right)\,.
\end{equation}
But $c_0(\epsilon)$ being $0$ at the ends of the interval $(0,  1-4/{\bf j}^2)) $, the maximum is obtained inside by looking at the zero of the derivative with respect to $\epsilon$.
We get
\begin{equation}
\epsilon _{max}= (1 -{\bf j}^{-2} )  -  \sqrt{(1- {\bf j}^{-2})^2  - (1- 4{\bf j}^{-2})      }   = (1-{\bf j}^{-2}) - {\bf j}^{-2} \sqrt{1 + 2 {\bf j}^2}\,.
\end{equation}
and 
\begin{equation}\label{forc0}
 c_0 =  2^{-6}{\bf j}^{-2}  \pi^3 \left( ({\bf j}^4 + 10 {\bf j}^2 - 2 )    -2  (2 {\bf j}^2+1)  \sqrt{1 + 2 {\bf j}^2} \right)  \,.
\end{equation}
Numerics with ${\bf j}$ replaced by its approximation gives
$
 c_0 \sim 0.014\,$.
 This is extremely small  and very far from from the conjectured value $2$ !\\
 \begin{remark}
 One can actually in \eqref{forc0}  replace ${\bf j}^2$ by $\frac{A(\Omega)}{\pi} \lim\inf \frac{\mathfrak L_k}{k}$. The constant ${\bf j}^2$ appears indeed only through \eqref{z2}. Because of the monotonicity of $c_0$
  as a function of ${\bf j}^2$ which results of the definition of $c_0$ as a sup. any improvement of a lower bound for $\frac{A(\Omega)}{\pi} \lim\inf \frac{\mathfrak L_k}{k}$ will lead to a corresponding improvment of $c_0$.
 \end{remark}
{\bf Acknowledgements.\\} I would like to thank T. Hoffmann-Ostenhof for former discussions on the subject, V. Felli and L. Abatangelo for inviting me to give a course on the subject in Milano (February 2015) and the Isaac Newton Institute where the final version of this note was completed, the author being Simons Foundation visiting  Fellow there.\\


\begin{thebibliography}{1}

 \bibitem  {AFT}
 B.~Alziary, J.~Fleckinger-Pell{\'e}, P.~Tak{\'a}{\v{c}}.
\newblock Eigenfunctions and {H}ardy inequalities for a magnetic
  {S}chr\"odinger operator in {$\mathbb R\sp 2$}.
\newblock  Math. Methods Appl. Sci. {\bf 26}(13), 1093--1136  (2003).

\bibitem  {BHV} V. Bonnaillie-No\"el, B. Helffer, and G. Vial.
\newblock Numerical simulations for nodal domains and spectral minimal
  partitions.
\newblock   ESAIM Control Optim.  Calc.Var. 
DOI:10.1051/cocv:2008074  (2008).

\bibitem {BBO} B. Bourdin, D. Bucur, and E. Oudet.
\newblock Optimal partitions for eigenvalues.
\newblock SIAM J. Sci. Comp. 31(6) 4100-4114  (2009).


\bibitem {CL} L.A. Caffarelli and F.H. Lin.
\newblock An optimal partition problem for eigenvalues.
\newblock  Journal of Scientific Computing {\bf 31}(1/2), DOI: 10.1007/s10915-006-9114.8 (2007).

\bibitem {CTV} M. Conti, S. Terracini, and G. Verzini.
\newblock  A variational problem for the spatial segregation of reaction-diffusion systems.
\newblock  Indiana Univ. Math. J. 54  (2005),  p.~779-815.

\bibitem{Cou} R. Courant.
\newblock  Ein allgemeiner Satz zur Theorie der Eigenfunktionen selbstadjungierter Differentialausdr\"ucke,
\newblock Nachr. Ges. G\"ottingen (1923), 81-84.

\bibitem{CH} R. Courant and D. Hilbert.
\newblock Methods of Mathematical Physics, Vol. 1.
\newblock New York (1953).

\bibitem {HHmag} B. Helffer and T. Hoffmann-Ostenhof.
\newblock On a magnetic characterization of spectral minimal partitions.
\newblock  JEMS 15, 2081--2092 (2013).

 \bibitem {HHO}
 B.~Helffer and  T.~Hoffmann-Ostenhof.
 \newblock A review on large k minimal spectral  k-partitions and  Pleijel's Theorem.
\newblock Proceedings of the congress in honour of  J. Ralston. Contemporary Mathematics, Vol. 640 (in press).


\bibitem {HHOT}
 B.~Helffer, T.~Hoffmann-Ostenhof, and S.~Terracini.
\newblock Nodal domains and spectral minimal partitions.
\newblock Ann. Inst. H. Poincar\'e Anal. Non
  Lin\'eaire {\bf 26}, 101--138 (2009).

\bibitem {HHOO}
B.~Helffer, M.~Hoffmann-Ostenhof, T.~Hoffmann-Ostenhof, and M.~P. Owen.
\newblock Nodal sets for groundstates of {S}chr\"odinger operators with zero
  magnetic field in non-simply connected domains.
\newblock {\em Comm. Math. Phys.} 202(3) (1999), p.~629-649.
  
    
\bibitem {Len} C. L\'ena.
\newblock Eigenvalues variations for Aharonov-Bohm operators.
\newblock J. Math. Phys. 56, 011502 (2015).

\bibitem {NT} B. Noris and S. Terracini.
\newblock Nodal sets of magnetic Schr\"odinger operators of
Aharonov-Bohm type and energy minimizing partitions.
\newblock Indiana Univ. Math. J., 59(4) :1361-1403, 2010.



  
\bibitem  {Pl} A.~Pleijel. 
\newblock Remarks on Courant's nodal theorem.
\newblock Comm. Pure. Appl. Math. {\bf 9},   543--550 (1956).


\end{thebibliography}
\end{document}